\documentclass{amsart}
\usepackage{amscd}
\numberwithin{equation}{section}
\def\arrowh#1{%
\begin{picture}(0,0)(0,0)
  \put(0,0){\line(1,0){2}}
  \put(0,0){\vector(1,0){1}}
  \put(1,-0.3){\makebox{$#1$}}
\end{picture}
}
\def\arrowv#1{%
\begin{picture}(0,0)(0,0)
  \put(0,0){\line(0,1){2}}
  \put(0,0){\vector(0,1){1}}
  \put(-0.3,1){\makebox{$#1$}}
\end{picture}
}
\def\Hdots{%
\begin{picture}(0,0)(0,0)
   \multiput(-0.2,0)(0.2,0){3}{\Dot}
\end{picture}
}
\def\Vdots{%
\begin{picture}(0,0)(0,0)
   \multiput(0,-0.2)(0,0.2){3}{\Dot}
\end{picture}
}
\def\Dot{\circle*{0.1}}
\def\DOT{\circle*{0.2}}

\def\Istructure{%
\begin{picture}(7,7)(-1,-1)
\multiput(0,0)(2,0){3}{\multiput(0,0)(0,2){3}{\DOT}}
\put(0,0){\arrowh{x}}
\put(0,0){\arrowv{y}}
\put(2,0){\arrowh{y}}
\put(2,0){\arrowv{x}}
\put(4,0){\arrowv{y}}
\put(5,0){\Hdots}
\put(0,2){\arrowh{y}}
\put(0,2){\arrowv{x}}
\put(2,2){\arrowh{x}}
\put(2,2){\arrowv{y}}
\put(4,2){\arrowv{x}}
\put(5,2){\Hdots}
\put(0,4){\arrowh{x}}
\put(2,4){\arrowh{y}}
\put(0,5){\Vdots}
\put(2,5){\Vdots}
\put(4,5){\Vdots}
\put(5,4){\Hdots}
\end{picture}
}
\let\cal\mathcal

\def\Uscr{{\cal U}}

\let\blb\mathbb

\def \ZZ{{\blb Z}}

\def \NN{{\blb N}}
\def \RR{{\blb R}}

\def\Id{\text{\upshape{id}}}

\def\Sym{\mathop{\text{\upshape{Sym}}}\nolimits}

\def\im{\operatorname {im}}

\def\ker{\operatorname {ker}}

\def\r{\rightarrow}

\newtheorem{lemma}{Lemma}[section]
\newtheorem{proposition}[lemma]{Proposition}
\newtheorem{theorem}[lemma]{Theorem}
\newtheorem{corollary}[lemma]{Corollary}
\newtheorem{condition}[lemma]{Condition}

\theoremstyle{definition}

\newtheorem{example}[lemma]{Example}

{

\newtheorem{case}{Case}

}

\theoremstyle{remark}

\newdimen\uboxsep \uboxsep=1ex
\def\uboxn#1{\vtop to 0pt{\hrule height 0pt depth 0pt\vskip\uboxsep
\hbox to 0pt{\hss #1\hss}\vss}}

\def\uboxs#1{\vbox to 0pt{\vss\hbox to 0pt{\hss #1\hss}
\vskip\uboxsep\hrule height 0pt depth 0pt}}

\title{Semigroups of $I$-type}
\author{Tatiana Gateva-Ivanova}
\address{Section of Algebra\\ Institute of Mathematics\\ Bulgarian
  Academy of Sciences\\ 113 Sofia\\ Bulgaria}
\email{tatiana@bgearn.bitnet}
\thanks{During the time the work on this paper was done, the 
 first author was partially supported by the
  J. W. Fulbright Exchange Program and by the Bulgarian Ministry of
  Education Grant MM-2:91}
\author{Michel Van den Bergh}
\address{Limburgs Universitair Centrum\\ Departement WNI\\ Universitaire
Campus\\ 3590 Die\-pen\-beek \\ Belgium}
\thanks{The second author is a senior researcher at the
NFWO}
\email{vdbergh@luc.ac.be, 
http://www.luc.ac.be/Research/Algebra/}

\begin{document}
\begin{abstract}
  Assume that $S$ is a semigroup generated by $\{x_1,\ldots,x_n\}$,
  and let $\Uscr$ be the multiplicative free commutative semigroup
  generated by $\{u_1,\ldots,u_n\}$. We say that $S$ is of
  \emph{$I$-typ}e if there is a bijection $v:\Uscr\r S$ such that for
  all $a\in\Uscr$, $\{v(u_1a),\ldots
  v(u_na)\}=\{x_1v(a),\ldots,x_nv(a)\}$. This condition
  appeared naturally in the work on Sklyanin algebras by John Tate and
  the second author.

In this paper we show that the condition for a semigroup to be of
$I$-type is related to various other mathematical notions found in the
literature. In particular we show that semigroups of $I$-type appear
in the study of the settheoretic solutions of the Yang-Baxter
equation, in the theory of Bieberbach groups and in the study of
certain skew binomial polynomial rings which were introduced by the
first author.
\end{abstract}
\maketitle
\section{Introduction}
In the sequel $k$ will be a field.
Our starting point for this paper are certain semi-groups which were
introduced in \cite{GI3}. Let $X=\{x_1,\ldots,x_n\}$ be a set of
generators. In \cite{GI3} the first author considers semigroups $S$ of the
form $\langle X;R\rangle$ where $R$ is a set of quadratic relations 
\[
R=\{x_jx_i=u_{ij}\mid i=1,\ldots,n; j=i+1,\ldots,n\}
\]
satisfying
{\def\thelemma{(*)}
\begin{condition}
\begin{enumerate}
\item
$u_{ij}=x_{i'}x_{j'}$, $i'<j'$, $i'<j$.
\item
As we vary $(i,j)$, every pair $(i',j')$ occurs exactly once.
\item
The overlaps $x_kx_jx_i$ for $k>j>i$ do not give rise to new relations
in $S$.
\end{enumerate}
\end{condition}
}
The motivation for (*) is developed in \cite{GI3}. Condition
(*1) says that the semigroup algebra $kS$ is a \emph{binomial skew
polynomial ring}, so the theory of (non-commutative) Gr\"obner bases  applies to it. Condition
(*3) says that as sets
\[
S=\{x_1^{a_1}\cdots x_n^{a_n}\mid (a_1,\ldots,a_n)\in\NN^n\}
\]
Furthermore it is shown in \cite[Thm II]{GI3} that   (*2) is
equivalent with $kS$ being noetherian (assuming (*1,3)). 

However conditions (*1,2,3) are also natural for intrinsic reasons.
There are exactly as many monomials $x_j x_i$ with $j>i$ as there are
monomials $x_{i'}x_{j'}$ with $i'<j'$. This provides the motivation
for imposing (*2). Furthermore, it follows from  \cite[Thm 3.16]{GI3}
that (*1,2,3) imply $j,j'>i,i'$ for the relations in $R$. Thus conditions
(*1,2,3) are actually symmetric, in the sense that if they are
satisfied by $S=\langle X;R\rangle$ then they are also satisfied by
$S^\circ$.

The purpose of this paper is to show that the semigroups defined in the
previous paragraphs are intimately connected with various other
mathematical notions which are currently of some interest. In particular
we show that they are  related to
\begin{enumerate}
\item
Set theoretic solutions of the Yang-Baxter equation \cite{Drinfeld1}.
\item
Bieberbach groups \cite{Charlap}.
\item
Rings of $I$-type \cite{TVdB}.
\end{enumerate}
We will now sketch these connections. We start by proving 
the following proposition.
\begin{theorem}
\label{intrth1}
Assume that  $R$ satisfies (*1,2,3). Define $r:X^2\r X^2$ as follows~:
$r$ is the identity on quadratic monomials and if
$(x_jx_i=x_{i'}x_{j'})\in R$ then $r(x_jx_i)=x_{i'}x_{j'}$,
$r(x_{i'}x_{j'})=x_jx_i$. 
Then $r$ satisfies
\begin{enumerate}
\item
$r^2=\Id_{X^2}$.
\item
$r$ satisfies the settheoretic Yang Baxter equation. That is, one has
\[
r_{1}r_{2}r_{1}=r_{2}r_{1}r_{2}
\]
where as usual 
$r_{i}:X^m\r X^m$ is defined as $\Id_{X^{i-1}}\times r\times \Id_{X^{m-i-1}}$.
\item
Given $a,b\in\{1,\ldots,n\}$ there exist unique $c,d$ such that
\[
r(x_cx_a)=x_dx_b
\]
Furthermore if $a=b$ then $c=d$.
\end{enumerate}
\end{theorem}
In view of this theorem it is natural to consider semigroups of the
form $\langle X;x_ix_j=r(x_ix_j)\rangle$ where $r$ is a
settheoretic solution of the Yang-Baxter equation. We will show that some
of these are of ``$I$-type'' \cite{TVdB}.
Being of $I$-type is a  technical condition which is very
useful for computations. Let us recall the definition here. We start
with a set of variables $u_1,\ldots,u_n$ and we let $\Uscr$ be the
free \emph{commutative} multiplicative semigroup generated by
$u_1,\ldots,u_n$. Let $S$ be a semigroup generated by $X=\{x_1,\ldots,x_n\}$. $S$ is said to be of
(left) $I$-type if there exists a bijection $v:\Uscr\r S$ (an
$I$-structure) such that $v(1)=1$ and such that for all $a\in\Uscr$
\begin{equation}
\label{intreq1}
\{v(u_1a),\ldots,v(u_na)\}=\{x_1 v(a),\ldots,x_nv(a)\}
\end{equation}
It is clear that if $S$ is of $I$-type then $kS$ is of $I$-type in the
sense of \cite{TVdB}. 

Assume that $S$ is $I$-type  with $I$-structure $v$. \eqref{intreq1}
implies that for every $a\in\Uscr$, $i\in\{1,\ldots,n\}$ there exists a 
unique $x_{a,i}\in X$ such that 
\[
x_{a,i}v(a)=v(au_i)
\]
and $\{x_{a,i}\mid i=1,\ldots, n\}=X$.
\begin{example}
\label{intrex14}
Let  $S$ be the semigroup $\langle x,y; x^2=y^2\rangle$ and consider the
 following doubly infinity graph.
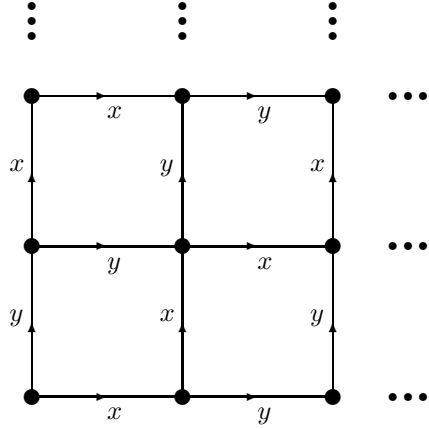
\begin{center}
\Istructure
\end{center} 
Define $v(u_1^{a_1} u_2^{a_2})$ as one (or all) of the paths 
from  $(0,0)$ to $(a_1,a_2)$, written in reverse order (for example
$v(u_1^2u_2)=xy^2=x^3=y^2x$).  Then it is clear that
this $v$ defines a $I$-structure on $S$.
\end{example}  We have the following result
\begin{theorem}
\label{intrth2} Assume that $S$ is $I$-type. Define
 $r:X^2\r X^2$ by 
\[
r(x_{u_i,j}x_{1,i})=x_{u_j,i}x_{1,j}
\]
Then $r$ satisfies the conclusions of Theorem \ref{intrth1}.
Conversely if $r:X^2\r X^2$ satisfies \ref{intrth1}.1.,2.,3. then the semigroup
$S=\langle X;x_i x_j=r(x_i x_j)\rangle$ is of $I$-type.
\end{theorem}
From Theorems \ref{intrth1},\ref{intrth2} it follows that semigroups defined by
relations satisfying (*1,2,3) are of $I$-type.
The proof of the following result is  similar
to the proof of \cite[Thm 1.1,1.2]{TVdB}.

For a cocycle $c:S^2\r k^\ast$ we use the notation $k_c S$ for the
twisted semi-group algebra associated to $(S,c)$. Thus $k_cS$ is the
$k$-algebra with basis $S$ and with multiplication $x\cdot
y=c(x,y)xy$ for $x,y\in S$.
\begin{theorem} \label{intrth3}
Assume that $S$ is of $I$-type and let 
 $A=k_cS$ for some cocycle $c:S^2\r k^\ast$. Then
\begin{enumerate}
\item
$A$ has finite global dimension.
\item
$A$ is Koszul.
\item
$A$ is noetherian.
\item
$A$ satisfies the Auslander condition.
\item
$A$ is Cohen-Macaulay.
\item If $c$ is trivial then $k_cS$ is finite over its center.
\end{enumerate}
\end{theorem}
For the definition of ``Cohen-Macaulay'' and the ``Auslander condition'' see
\cite{Le}. 
\begin{corollary}
\label{intrcor}
Assume that $S$ is a semigroup of $I$-type. Then $k_cS$ is a domain, and in
particular $S$ is a cancellative. 
\end{corollary}
This corollary follows from \cite{Le}.

Let $S$ be a semi-group of $I$-type with $I$-structure $v:\Uscr\r S$.
Since $S$ is a cancellative semigroup of subexponential growth, it is
\"Ore. Denote its quotient group by $\bar{S}$.  We identify $\Uscr$
in the natural way with $\NN^n$, and in this way we  embed it in $\RR^n$. We
will prove the following
\begin{theorem}
\label{intrth4} Assume that $S$ is of $I$-type with $I$-structure
$v:\Uscr\r S$.  Let $S$ act on the right of $\Uscr$ by pulling back
under $v$ the action of $S$ on itself by right translation. Then this
action extends to a free right action of $\bar{S}$ on $\RR^n$ by Euclidean
transformations and for this action $[0,1[^n$ is a fundamental domain.
In particular $\bar{S}$ is a Bieberbach group.
\end{theorem}
\begin{example} If we take for $S$ the semigroup of Example \ref
{intrex14} then using \eqref{euaction} one checks that $x$ and $y$ act  
on $\RR^2$
 by
glide reflections along parallex axes. Hence $\RR^2/\bar{S}$ is the Klein
bottle!
\end{example}
\section{Proof of Theorem \ref{intrth1}}
In this section we prove Theorem \ref{intrth1}. The notations will be
as in the introduction. So $S$ is a semigroup of the form $\langle
X;R\rangle$ where $R$ is a set of relations satisfying~(*). It
is clear that \ref{intrth1}.1. is true by definition.
So we concentrate
on \ref{intrth1}.2. and \ref{intrth1}.3.

Below we denote the diagonal of $X^m$ by $\Delta_m$. Clearly
\[
r_1(\Delta_3)=\Delta_3,\qquad r_2(\Delta_3)=\Delta_3
\]
Furthermore it follows from the ``cyclic condition''
\cite[Thm 3.16]{GI3} that
\begin{equation}
\label{cyclic}
r_1r_2(\Delta_2\times X)=X\times \Delta_2
\end{equation}
\begin{lemma}
\label{somelemma}
The relation
\[
r(zt)=xy
\]
defines bijections between $X^2$ and itself given by
\[
 (t,y)\leftrightarrow(z,t)\leftrightarrow (x,y)\leftrightarrow (z,x)
\]
\end{lemma}
\begin{proof}
That $(z,t)\leftrightarrow (x,y)$ defines a bijection is clear. Now
consider the map which assigns $(t,y)$ to $(z,t)$. We claim that it is
an injection. If this is so then by looking at the cardinality of the
source and the target (which are both $X^2$) we see that it must be a
bijection. 

To prove the claim we compute $r_2r_1(xy^2)=r_2(zty)=z^2\ast$ where
the last equality follows from \eqref{cyclic}. Thus $r(ty)=z\ast$ and
hence $z$ is uniquely determined by $t,y$. This proves the claim.

That $(z,t)\leftrightarrow (z,x)$ is a bijection is proved similarly.
\end{proof}
Note that lemma \ref{somelemma} contains  \ref{intrth1}.3 as a special case.
Hence we are  left with proving \ref{intrth1}.2.

Let us call $w,w'\in \langle X\rangle $ equivalent if they have the
same image in  $S$. Notation~: $w\sim w'$.  Clearly $w\sim w'$ iff
\[
w'=r_{i_1}r_{i_2}\cdots r_{i_p} w
\]
for some $p,i_1,\ldots,i_p$.

Concerning the structure of the
equivalence classes there is the following easy lemma.
\begin{lemma}
\label{easy}
Every equivalence class for $\sim$ in $X^m$ contains exactly one
monomial of the form $x_{a_1}\cdots x_{a_m}$, $a_1\le \cdots \le a_m$.
\end{lemma}
\begin{proof}
This is a consequence of the Bergman diamond
lemma.
\end{proof}
After these preliminaries we  prove the Yang-Baxter relation for
$r$. The proof is based upon a careful examination of  the
equivalence classes in $X^3$, together with a counting argument.

Let $D$ be the infinite dihedral group
$\langle r_1,r_2; r_1^2=r_2^2=e\rangle$. $D$ acts on $X^3$ and is clear
the the equivalence classes correspond to $D$-orbits. Let $O$ be such
an orbit. There are three possibilities.
\begin{itemize}
\item[(A)] $O\cap\Delta_3\neq \emptyset$. In this case clearly $|O|=1$.
\item[(B)] $O\cap ((\Delta_2\times X\cup X\times \Delta_2)\setminus
  \Delta_3)\neq \emptyset$.
  In this case it follows from \eqref{cyclic} that $|O|=3$.
\item[(C)] $O\cap (\Delta_2\times X\cup X\times\Delta_2)=\emptyset$.
Now $O=\{w,r_1w,r_2r_1w,\ldots\}$. Thus a general member of $O$ is of the
form $(r_2r_1)^aw$ or $r_1(r_2r_1)^a w$. 

We claim that  $(r_2r_1)^aw\neq r_1(r_2r_1)^b w$ for $a,b\in\ZZ$. To
prove this, assume the contrary and
 define \[
w_1=\begin{cases}
r_1(r_2r_1)^{\left\lfloor\frac{a+b}{2}\right\rfloor}w&\text{if $a+b$ is odd}\\
(r_2r_1)^{\left\lfloor\frac{a+b}{2}\right\rfloor}w&\text{if $a+b$ is even}
\end{cases}
\]
Thus $r_1w_1=w_1$ or $r_2w_1=w_1$ (depending on whether $a+b$ is even
or odd), whence $w_1\in \Delta_2\times X\cup X\times \Delta_2$,
contradicting the hypotheses.

Let $p$ be the smallest positive integer such that $(r_2r_1)^pw=w$. Then
\begin{equation}
\label{ybeq1}
O=\{w,(r_2r_1)w,\ldots,(r_2r_1)^{p-1}w,r_1w,r_1(r_2r_1)w,\ldots,r_1
(r_2r_1)^{p-1}w\}
\end{equation}
In particular $|O|=2p$ is even. We claim $|O|\ge 6$. To prove this we have
to exclude $|O|=2,4$. The case $|O|=2$ is easily excluded using
\ref{intrth1}.3. Hence we are left with $|O|=4$. This means that $O$ looks
like
\[
\begin{CD}
x_ax_bx_c @>r_2 >> x_a x_d x_e\\
@V r_1VV @VV r_1 V\\
x_f x_g x_c @>r_2 >> x_f x_h x_e
\end{CD}
\]
which implies that $R$ contains relations
\begin{align}
\label{al1}x_b x_c&=x_d x_e\\
\label{al2}x_ax_b&=x_fx_g\\
\label{al3}x_a x_d&=x_fx_h\\
\label{al4}x_g x_c&=x_hx_e
\end{align}
Now in a relation $x_ux_v=x_wx_t$ the couples $(u,v)$ and $(v,t)$
determine each other (lemma \ref{somelemma}). So looking at
\eqref{al2}\eqref{al3} we find $b=d$, $g=h$. This implies that
\eqref{al1} is actually of the form $x_dx_c=x_dx_c$, which is a
contradiction. Hence $|O|\ge 6$.
\end{itemize}
An alternative classification of these orbits goes through the
elements they contain of the form $x_ax_b x_c$, $a\le b\le c$. A
unique such element exist in every orbit by lemma~\ref{easy}.  

If
 $O$ contains an element of the form $x_ax_bx_c$, $a<b<c$ then it
is of type (C) because if not,  it contains an element of the
form $x_dx_dx_e$ or $x_dx_ex_e$ with $d\ge e$. Using \eqref{cyclic} and
(*1) such elements are equivalent to elements of the form $x_fx_g
x_g$, $x_fx_fx_g$ with $f\le g$. Contradiction.

If $O$ contains an element of the form $x_ax_ax_b$ or of
the form $x_ax_bx_b$ with $a<b$ then $O$ is clearly of type
(B). Finally $O$ is of type (A) iff it contains an element of the form
$x_ax_ax_a$. 

Thus we find that there are $n$ orbits of type (A), $n(n-1)$ orbits of
type (B) and $n(n-1)(n-2)/6$ orbits of type (C). From the equality
\[
|X^3|=n^3=1\cdot n+2\cdot n(n-1)+6\cdot \frac{n(n-1)(n-2)}{6}
\]
we deduce that the orbits of type (C) contain exactly $6$ elements.

Now Yang-Baxter easily follows. If $w$ has orbit of type (C) then
from \eqref{ybeq1} we deduce that $(r_2r_1)^3w=w$. If the orbit is of
type (B) then $(r_2r_1)^3w=w$ follows directly from
\eqref{cyclic}. Finally if the orbit is of type (A) then $r_1w=r_2w=w$
and there is nothing to prove.

This concludes the proof of Theorem \ref{intrth1}.
\section{Proof of Theorem \ref{intrth2}}
In this section we prove Theorem \ref{intrth2}. One direction is
trivial, so we concentrate on the other direction. That is, given $r$
satisfying \ref{intrth1}.1.,2.,3. we will construct $v:\Uscr\r S$
and $x_{b,i}\in X$ for $b\in\Uscr$, $i=\{1,\ldots,n\}$ in such a way
that
\begin{itemize}
\item[(a)] $v$ is a bijection.
\item[(b)] $v(u_i b)=x_{b,i}v(b)$
\item[(c)] $\{x_{b,i}\mid i=1,\ldots,n\}=\{x_1,\ldots,x_n\}$
\item[(d)] $r(x_{bu_j,i}x_{b,j})=x_{bu_i,j}x_{b,i}$
\end{itemize}
The construction is inductive. To start we put $v(1)=1$ and
$v(u_i)=x_{\sigma(i)}$ for an  arbitrary element $\sigma$ of $\Sym_n$.
From here on everything will be uniquely defined.
Assume  that we have constructed $v(b)$ for $\deg b\le m-1$, 
$x_{b,i}$ for $\deg b\le m-2$ satisfying (a-d). We will define $x_{a,i}$
for $\deg a=m-1$ such that (c)(d) hold.
\begin{case} $a\neq u_i^{m-1}$. So $a=bu_j$, $j\neq i$. Computing
  $v(bu_iu_j)$ in two ways (as a heuristic device, since $v(bu_iu_j)$
  is still undefined) we find that $x_{a,i}$ must be defined by
\begin{equation}
\label{alpha}
r(x_{a,i}x_{b,j})=\ast x_{b,i}
\end{equation}
This indeed defines $x_{a,i}$  uniquely thanks to
\ref{intrth1}.3. However one still must deal with the possibility that
$x_{a,i}$ might depend on $j$. To analyze this assume $k\neq i$,
$a=du_ju_k$. Put $b=du_k$, $c=du_j$, $e=du_i$. We now define
$p,q,p',q'$ by
\begin{align}
\label{sec2eq1} r(px_{b,j})&=qx_{b,i}\\
\label{sec2eq2}
r(p'x_{c,k})&=q'x_{c,i}
\end{align}
We have to show $p=p'$. By induction we have the following identities.
\begin{align}
r(x_{b,j}x_{d,k})&=x_{c,k}x_{d,j}\\
r(x_{b,i}x_{d,k})&=x_{e,k}x_{d,i}\\
\label{comp} r(x_{c,i}x_{d,j})&= x_{e,j}x_{d,i}
\end{align} 
We can now construct a ``Yang-Baxter diagram''
\[
\begin{CD}
px_{b,j}x_{d,k} @>r_1>> qx_{b,i} x_{d,k}\\
@V r_2VV @VV r_2 V\\
px_{c,k} x_{d,j} @. qx_{e,k} x_{d,i}\\
@V r_1 VV @VV r_1V\\
X Y x_{d,j} @>r_2 >> X Z x_{d,i}
\end{CD}
\]
with $X,Y,Z$ unknown sofar.

Comparing $r(Yx_{d,j})=Zx_{d,i}$ with \eqref{comp} yields
$Y=x_{c,i}$, $Z=x_{e,j}$. 

So we find that
\[
r(px_{c,k})=Xx_{c,i}
\]
and comparing this with \eqref{sec2eq2} yields $p=p'$.

Hence we can now legally define $x_{a,i}=p$. Furthermore
\eqref{sec2eq1} can also be read as
\[
r(qx_{b,i})=px_{b,j}
\]
Since obviously $bu_i\neq u^{m-1}_j$ we obtain $q=x_{{bu_i},j}$. We
conclude that with our present definitions we have for $j\neq i$, $\deg b\le m-2$
\begin{equation}
\label{sec2eq3}
r(x_{bu_j,i}x_{b,j})=x_{bu_i,j}x_{b,i}
\end{equation}
We claim that this relation  holds more generally under the
hypotheses that $\deg b\le m-2$ and $bu_j\neq u_i^{m-1}$ (or equivalently
 $bu_i\neq u_j^{m-1}$).

The only case that still has to be checked is~: $i=j$, $\deg b=m-2$,
$b\neq u_i^{m-2}$.
In this case we may put $b=cu_k$, $k\neq i$. We construct
again a Yang-Baxter diagram
\begin{equation}
\label{beta}
\begin{CD}
x_{cu_iu_k,i} x_{cu_k,i} x_{c,k} @>r_1>> x_{cu_iu_k,i} Y x_{c,k}\\
@V r_2 VV @AA r_2 A\\
x_{cu_iu_k,i}x_{cu_i,k}x_{c,i} @. x_{cu_iu_k,i}x_{cu_i,k} x_{c,i}\\
@V r_1 VV  @AA r_1A\\
x_{cu_i^2,k} x_{cu_i,i} x_{c,i} @>r_2>> x_{cu^2_i,k} x_{cu_i,i} x_{c,i}
\end{CD}
\end{equation}
From the relation 
\[
r(x_{cu_i,k}x_{c,i})=Yx_{c,k}
\]
we deduce $Y=x_{cu_k,i}$. 
Looking at the toprow of \eqref{beta} finishes the proof
of \eqref{sec2eq3} under the hypotheses that $bu_j\neq u_i^{m-1}$.

Now we claim that if $\deg a=m-1$, $i\neq j$ and $a\neq u_i^{m-1},
u_j^{m-1}$ then $x_{a,i}\neq x_{a,j}$. Assume the contrary and write 
$a=bu_l$.
Then by \eqref{sec2eq3} we have
\begin{equation}
\label{sec2eq4}
\begin{split}
r(x_{bu_l,i}x_{b,l})&=x_{bu_i,l} x_{b,i}\\
r(x_{bu_l,j}x_{b,l})&=x_{bu_j,l}x_{b,j}\\
\end{split}
\end{equation}
Since the lefthand sides of \eqref{sec2eq4} are the same and this is
not the case with the righthand sides we obtain a
contradiction.
\end{case}
\begin{case}
$a=u_i^{m-1}$. In this case we take $x_{a,i}$ different from $x_{a,j}$,
$j\neq i$. This defines $x_{a,i}$ uniquely, and obviously (c) is
satisfied if $\deg b\le m-1$.

Now we prove \eqref{sec2eq3} in the remaining case $b=u_i^{m-2}$, $i=j$.

Since we already know (c) we can
write
\[
r(x_{bu_k,l}x_{b,k})=x_{bu_i,i}x_{b,i}
\]
for some $k$, $l$ and we have to show $k=l=i$. Assume on the contrary
that $k\neq i$ or $l\neq i$. By what we know so far we have
\[
r(x_{bu_k,l}x_{b,k})=x_{bu_l,k}x_{b,l}
\]
But then $k=l=i$. Contradiction.
\end{case}
So up to this point we have defined $x_{b,i}$ and we have proved (c)(d)
for $\deg b\le m-1$. Now if $a=bu_i$ has length $m$ then we define
\begin{equation}
\label{sec2eq5}
v(a)=x_{b,i}v(b)
\end{equation}
so that (b) certainly holds. That \eqref{sec2eq5} is well defined follows
easily from (d). 

Hence to complete the induction step it suffices to show that (a)
holds. That is $v$ should define a bijection on words of length $m$. 
Let $U=\{u_1,\ldots,u_n\}$ and let $U^m$ be the words of length  $m$ in
$U$. Furthemore let $r_i:U^m\r U^m$ be given by exchanging the $i,i+1$'th
letter. Define a map
$
\tilde{v}:U^m\r X^m
$
by
\[
\tilde{v}(u_{i_1}\cdots u_{i_m})=x_{u_{i_2}\cdots u_{i_m},i_1}\cdots
x_{u_{i_{m-1}}u_{i_m},i_{m-2}}x_{u_{i_m},i_{m-1}}x_{1,i_m}
\]
By (c), $\tilde{v}$ is clearly a bijection.

From (d) we obtain the following commutative diagram.
\[
\begin{CD}
U^m @>\tilde{v}>> X^m\\
@V r_{i} VV @VV r_{i} V\\
U^m @>\tilde{v} >> X^m
\end{CD}
\]
So $\tilde{v}$ defines a bijection between the orbits $U^m/\Sym_m$ and
$X^m/\Sym_m$. We have
\[
\Uscr_m=U^m/\Sym_m,\qquad S_m=X^m/\Sym_m
\]
where $\Uscr_m$, $S_m$ are the elements of degree $m$ in $\Uscr$ and $S$
respectively. Furthermore the map $\Uscr_m\r S_m$ induced by $\tilde{v}$
is precisely $v$. This finishes the proof of Theorem~\ref{intrth2}.

\section{Semigroups of $I$-type}
Below $S$ will be a semigroup of $I$-type, with $I$-structure
$v:\Uscr\r S$ (as defined in the  introduction). In this section we will
give some properties of $S$, and in particular we will prove Theorem
\ref{intrth3}.

First observe that every element of $\langle X\rangle$ can be written
uniquely in the form 
\[
x_{u_1\cdots u_{i_{m-1}},i_{m}}\cdots
x_{u_{i_1},i_2} x_{1,i_1}
\]
Two different elements $w$, $w'$ in $X^2$ have the same image in $S$ iff
there exist $i\neq j$ such that
\[
w=x_{u_i,j}x_{1,i},\qquad w'=x_{u_j,i}x_{1,j}
\]
The following lemma summarizes some observations in \cite{TVdB}, translated
into the language of semigroups.
\begin{lemma}
\label{ZZlemma41}
\begin{enumerate}
\item The natural grading by degree on $\Uscr$ induces via $v$ a grading
on $S$ such that $\deg(x_i)=1$.
\item The map $s\mapsto sv(\mu)$ for a given $\mu\in\Uscr$ induces a
bijection between $S$ and $\{v(a\mu)\mid a\in\Uscr\}$.
\item $S$ is right cancellative.
\item $S$ is a quotient of $\langle X\rangle$ by $n/(n-1)/2$ different
relations in degree $2$ given by
\[
x_{u_i,j}x_{1,i}=x_{u_j,i}x_{1,j},\qquad j>i
\]
\end{enumerate}
\end{lemma}

\noindent
If $\sigma\in \Sym_n$ then we extend $\sigma$ to $\Uscr$ via 
\[
\sigma (u_{i_1}\cdots u_{i_p})=u_{\sigma{i_1}}\cdots u_{\sigma i_p}
\]
\begin{lemma}
\label{ZZlemma42}
Every bijection $w:\Uscr\r S$, satisfying \eqref{intreq1} is of the form
$v\circ \sigma$, $\sigma\in \Sym_n$.
\end{lemma}
\begin{proof} Clearly there exist $\sigma \in \Sym_n$ such that $w$ and
$v\circ \sigma$ take the same values on $\{u_1,\ldots,u_n\}$. Hence to
prove the lemma we have to show that a map $v$ satisfying
\eqref{intreq1} is uniquely determined by the values it takes on
$\{u_1,\ldots,u_n\}$. This  was part of the proof of Theorem
\ref{intrth2}.
\end{proof}
Now we want to develop some kind of calculus for semigroups of $I$-type.
Consider the arrows
\begin{equation}
\label{somearrows}
\begin{CD}
S @>s\mapsto sv(b)>> \{v(ab)\mid b\in\Uscr\} \\
@. @AA \begin{matrix} v(ab)\\ \uparrow\\b\end{matrix} A\\
@. \Uscr
\end{CD}
\end{equation}
It is clear that the vertical map is a bijection and
so is the horizontal map by lemma \ref{ZZlemma41}.
Thus  we may define a
bijection $w:\Uscr\r S$ which makes \eqref{somearrows}
commutative. Furthermore $w$ obviously satisfies \eqref{intreq1}, so
according to lemma \ref{ZZlemma42} $w=v\circ \phi(b)$ where
$\phi(b)\in\Sym_n$. We view $\phi$ as a map from $\Uscr$ to $\Sym_n$.
Expressing the fact that $w$ completes \eqref{somearrows} to a commutative diagram yields
\begin{equation}
\label{basiceq1}
v(ab)=v(\phi(b)(a))\,v(b)
\end{equation}
If we now compute $v(abc)$ in two ways we find
\[
v(abc)=v(\phi(\phi(c)(b))(\phi(c)(a)))\,v(\phi(c)(b))\,v(c)
\]
and
\[
v(abc)=v(\phi(bc)(a))\,v(\phi(c)(b))\,v(c)
\]
Using the fact that $S$ is right cancellative we obtain
\[
\phi(\phi(c)(b))(\phi(c)(a))=\phi(bc)(a)
\]
or put differently
\[
(\phi(\phi(c)(b))\circ \phi(c))(a)=\phi(bc)(a)
\]
Since this is true for all $a$ be obtain 
 \begin{equation}
\label{basiceq2}
\phi(bc)=\phi(\phi(c)(b))\circ \phi(c)
\end{equation}
Let us define $\ker\phi$, $\im \phi$ in the usual way (even though
$\phi$ is clearly not a semigroup homomorphism).
\begin{align*}
\ker\phi &=\{a\in\Uscr\mid \phi(a)=\Id\}\\
\im\phi&=\{\phi(a)\mid a\in\Uscr\}
\end{align*}
To simplify the notation we put $P=\ker\phi$, $G=\im\phi$.

Then \eqref{basiceq1}\eqref{basiceq2} yield the following lemma.
\begin{lemma}
\begin{enumerate}
\item If $b\in P$ then 
\begin{align}
\label{eqa}\phi(ab)&=\phi(a)\\
\label{eqb} v(ab)&=v(a)v(b)
\end{align}
\item $P$ is a saturated subsemigroup of $\Uscr$ ($a\in P\Rightarrow
  (ab\in P\iff b\in P)$).
\item $G$ is a subgroup of $\Sym_n$ (note that a finite subsemigroup 
of a group is itself a group).
\item If $b\in G$ and $a\in P$ then $b(a)\in P$. 
\end{enumerate}
\end{lemma}
\begin{lemma} 
\label{ZZlem43}
There exist $t_1,\ldots,t_n>0$ such that
  $u_i^{t_i}\in P$.
\end{lemma}
\begin{proof} Since $\Sym_n$ is finite there exist $r_i<s_i$ such that
\begin{equation}
\label{secXXeq1}
\phi(u_i^{r_i})=\phi(u_i^{s_i})
\end{equation}
Put $a=\prod_i u_i^{r_i}$, $t'_i=s_i-r_i$.

Now if $\phi(p)=\phi(q)$ then \eqref{basiceq2} implies that
$\phi(rp)=\phi(rq)$. Applying this with $p=u^{r_i}$, $q=u^{s_i}$ and
$r=\prod_{j\neq i} u_j^{r_i}$ yields
$\phi(a)=\phi(au^{t'_i})=\phi(\phi(a)(u^{t_i}_i))\phi(a)$ and thus 
\[
\phi(a)(u_i)^{t'_i}\in\ker\phi
\]
Now $\phi(a)(u_i)=u_{\phi(a)(i)}$ so if we put $t_i=t'_{\phi(a)(i)}$ then
$\phi(u_i^{t_i})=\Id$.
\end{proof}
\begin{corollary}
\label{eencor}
Let $P_0$ be the subsemigroup of $\Uscr$ generated by $u^{t_i}_i$. Then 
\begin{enumerate}
\item $v(P_0)$ is a free abelian subsemigroup of $S$, generated by
$v(u_i^{t_i})$.
\item
$S=\bigcup_a v(a) v(P_0)$
\end{enumerate}
where the union runs over those $a=u_1^{p_1}\cdots u_n^{p_n}$ with $0\le
p_i\le t_i-1$.
\end{corollary}
\begin{proof}
The corresponding statements for $\Uscr$ are obvious. To obtain them for
$S$ one
applies $v$ and  uses \eqref{eqb}. \end{proof}
\begin{proof}[Proof of Theorem \ref{intrth3}] This is entirely similar to
  the proof of \cite[Thm 1.1, Thm 1.2]{TVdB} so we content ourselves
  with a quick sketch.  Note that by \cite[Cor 3.6]{TVdB} an algebra of
  $I$-type is automatically Koszul and has finite global dimension, so
  we only have to prove 3.-6.

  Note that the equations of $k_cS$ are given by
  $x_{u_i,j}x_{1,i}=d_{ij}x_{u_j,i}x_{1,j}$ for some $d_{ij} \in
  k^\ast$.  We first assume that the $(d_{ij})_{ij}$ are roots of
  unity. Then (using \eqref{eqb}) we can take $P_0$ so small that
  $v(P_0)$ is commutative in $k_cS$. Thus by corollary \ref{eencor},
  $k_cS$ is finite on the left over a commutative ring, and hence is
  PI. This proves in particular 6.  and using the same results of
  Stafford and Zhang \cite{Staf2} as in the proof of \cite[Thm
  1.1]{TVdB} also yields 2.-5. in this case.

The general case is now proved using reduction to a finite field as in
\cite{TVdB}.
\end{proof}
\section{Proof of Theorem \ref{intrth4}}
In this section we use the same notations and assumptions as in the
previous sections.

Since $S$ is cancellative (Cor. \ref{intrcor}) and has subexponential
growth it is (left and right) \"Ore.
For an \"Ore semigroup $T$ denote by $\bar{T}$ its quotient group.

We now extend $v$, $\phi$ to maps
\begin{align*}
\bar{v}&:\bar{\Uscr}\r \bar{S}:up^{-1}\mapsto v(u)v(p)^{-1}\\
\bar{\phi}&:\bar{\Uscr}\r \Sym_n:up^{-1}\mapsto \phi(u)\phi(p)^{-1}
\end{align*}
where $p\in P$. This is well defined because of \eqref{eqa}\eqref{eqb} and
the fact  that it is clear from lemma \ref{ZZlem43} that every element of
$\bar{\Uscr}$ can be written as $up^{-1}$, $p\in P_0\subset P$.
\begin{lemma}
\begin{enumerate}
\item
If $s\in S$ then there exists $t\in S$ such that $ts\in v(P)$, $st\in v(
P)$. \item
$\bar{v}$ is a bijection.
\end{enumerate}
\end{lemma}
\begin{proof}
\begin{enumerate}
\item
Assume $t=v(c)$. We have to find $b\in\Uscr$ such that 
\begin{align*}
\phi(v^{-1}(v(b)v(c)))=\phi(b)\phi(c)=\Id\\
\phi(v^{-1}(v(c)v(b)))=\phi(c)\phi(b)=\Id
\end{align*}
It is clear that this is possible since $\im\phi$ is a group.
\item
It is easy to see that $\bar{v}$ is an injection, and from 1.\ we
deduce that it is also a surjection. \qed \end{enumerate}
\def\qed{}\end{proof}
One verifies that $\bar{v}$ satisfies \eqref{intreq1} and it is also
clear $\ker\bar{\phi}$, $\im\bar{\phi}$ have the same properties as
$\ker\phi$, $\im \phi$ (lemma \ref{ZZlem43}). Furthermore $\ker
\bar{\phi}$ is now actually a group and $\im\bar{\phi}=\im\phi$. We
deduce the following slight strengthening of lemma \ref{ZZlem43} (and
generalization of \cite{GI3}) which is however not needed in the
sequel.
\begin{proposition} For all $i$~: $u_i^{n!}\in \ker\phi$.
\end{proposition}
\begin{proof} Let $p$ be the smallest positive integer such that
$u^p_i\in \ker\phi
$. Then 
$
p$ divides $| \bar{\Uscr}/\ker \bar{\phi}|
$
Now $\bar{\phi}$ defines a bijection (\emph{not} a group homomorphism)
between $\bar{\Uscr}/\ker\bar{\phi}$ and $\im\bar{\phi}$.
Thus
$
p$ divides $|\im \bar{\phi}|$ which in turn divides $|\Sym_n|=n!$
\end{proof}
$\bar{S}$ acts on itself by right and left multiplictation. If we
transport this action to $\Uscr$ through $v$ we obtain \emph{commuting}
left and right actions of $\bar{S}$ on $\bar{\Uscr}$ given by the formulas
\begin{align}
\label{action1}
\forall a\in\bar{S}, b\in\bar{\Uscr}&:a\cdot
b=\bar{v}^{-1}(a\bar{v}(b))\\
\label{action2}
\forall a\in\bar{\Uscr},b\in\bar{S}&:a\cdot b=\bar{v}^{-1}(\bar{v}
(a) b) \end{align}
In the previous sections we have concentrated on the action
\eqref{action1}. Now we will say something about the action \eqref{action2}.

Using \eqref{basiceq1} we deduce that for $a\in\bar{\Uscr}$, $b\in
\bar{S}$~: \[
a\cdot b= \bar{\phi}(\bar{v}^{-1}(b))^{-1}(a)\, \bar{v}^{-1}(b)
\]
\begin{proof}[Proof of Theorem \ref{intrth4}] By  permuting the $x_i$
we may and we will assume that $v(u_i)=x_i$.
Consider the map
\[
\psi:\ZZ^n\r \bar{\Uscr}:(a_1,\ldots,a_n)\mapsto u_1^{a_1}\cdots
u_n^{a_n}
\]
For $a\in\ZZ^n$, $b\in\bar{S}$ we write
\[
a\cdot b=\psi^{-1}(\psi(a)\cdot b)
\]
and we put $\tilde{\phi}(c)=\phi(c)\circ \psi$,
$\tilde{\phi}_i=\tilde{\phi}(u_i)$. We find for $(a_1,\ldots,a_n)\in\ZZ^n$~:
\begin{equation}
\label{euaction}
(a_1,\ldots,a_n)\cdot x_i=
(a_{\tilde{\phi}_i(1)},\ldots, a_{\tilde{\phi}_i(i)}+1,\ldots,
a_{\tilde{\phi}_i(n)})
\end{equation}
We conclude that $(x_i)_i$, and hence all of $\bar{S}$ acts on the 
right of $\ZZ^n$ by Euclidean transformations. Keeping the formula
\eqref{euaction} we can extend this action to an action on $\RR^n$ and
it is then clear  that $[0,1[^n$ is a fundamental domain. Furthermore if
the action were not free then there would be a fixed point
$(a_1,\ldots,a_n)\in\RR^n$ for some element $s$ of $\bar{S}$. But then
$(\lfloor a_1\rfloor,\ldots,\lfloor a_n\rfloor )\in\ZZ^n$ is also a fixed
point for $s$.  This is impossible since by construction the action of $\bar{S} $ on $\Uscr$
and hence on $\ZZ^n$ is free.
\end{proof}
\ifx\undefined\bysame
\newcommand{\bysame}{\leavevmode\hbox to3em{\hrulefill}\,}
\fi

\end{document}